\documentclass[a4paper,12pt]{amsart}

\usepackage{amstext}
\usepackage{amsmath}
\usepackage{ifthen}
\usepackage{amscd}
\usepackage{amssymb}
\usepackage[french]{babel}
\usepackage[T1]{fontenc}
\usepackage[utf8]{inputenc}
\usepackage[all]{xy}
\usepackage{graphicx}
\usepackage{enumerate}
\usepackage{xspace}
\usepackage{pdfsync}
\usepackage{epic}
\usepackage{dsfont}
\usepackage{stackrel}

\newenvironment{demo}[1][]{\ifthenelse{\equal{#1}{}}{\noindent\textbf{Démonstration :}\xspace}{\noindent\textbf{Démonstration #1:}\xspace}}{$\square$\newline}

\newtheoremstyle{break}
  {}
  {}
  {\itshape}
  {}
  {\bfseries}
  {.}
  {\newline}
  {}
  
\newtheoremstyle{rq}
  {}
  {}
  {\slshape}
  {}
  {\bfseries}
  {.}
  {3pt}
  {}

\newtheoremstyle{exemple}
  {}
  {}
  {\upshape}
  {}
  {\bfseries}
  {.}
  {3pt}
  {}

\newtheoremstyle{fact}
  {}
  {}
  {\slshape}
  {}
  {\bfseries}
  {.}
  {2pt}
  {}

\theoremstyle{fact}
\newtheorem*{fact*}{Fact}
\theoremstyle{break}

\newtheorem{thm}{Théorème}[section]

\newtheorem{cor}[thm]{Corollaire}
\newtheorem{lem}[thm]{Lemme}
\newtheorem{prop}[thm]{Proposition}
\newtheorem{defi}[thm]{Définition}
\theoremstyle{rq}
\newtheorem{rem}[thm]{Remarque}
\newtheorem{qt}[thm]{Question}
\theoremstyle{exemple}

\newtheorem*{merci}{Remerciements}

\newcommand{\RR}{\mathbb R}

\newcommand{\CC}{\mathbb C}
\newcommand{\QQ}{\mathbb Q}
\newcommand{\ZZ}{\mathbb Z}

\newcommand{\PP}{\mathbb P}

\newcommand{\To}{\longrightarrow}

\title[Groupes kählériens et projectifs]{Repr\'esentations lin\'eaires des groupes k\"ahlériens et de leurs analogues projectifs}
\date{\today}
\author{Fréderic Campana, Beno\^it Claudon, Philippe Eyssidieux}
\address{Fréderic \textsc{Campana}, Université de Lorraine, Institut \'Elie Cartan Nancy, UMR 7502, B.P. 70239, 54506 Vand\oe uvre-lès-Nancy Cedex, France}
\address{Beno\^it \textsc{Claudon}, UMI CNRS/IMPA, Estrada Dona Castorina 110, Jardim Botânico, 22460-320, Rio de Janeiro, Brasil}
\address{Philippe \textsc{Eyssidieux}, Institut Fourier, Université Grenoble 1, 38402 Saint-Martin d'Hères Cedex, France}
\email{Frederic.Campana@univ-lorraine.fr}
\email{Benoit.Claudon@univ-lorraine.fr}
\email{Philippe.Eyssidieux@ujf-grenoble.fr}

\begin{document}

\begin{abstract}
Dans cette note nous \'etablissons le r\'esultat suivant, annonc\'e dans \cite{CCE} : si $G\subset \mathrm{GL}_n(\Bbb C)$ est l'image d'une repr\'esentation lin\'eaire d'un groupe k\"ahlérien $\pi_1(X)$, il admet un sous-groupe d'indice fini qui est l'image d'une repr\'esentation lin\'eaire du groupe fondamental d'une vari\'et\'e projective complexe lisse $X'$.

Il s'agit donc de la solution (\`a indice fini pr\`es) pour les repr\'esentations lin\'eaires d'une question usuelle demandant si le groupe fondamental d'une vari\'et\'e k\"ahl\'erienne compacte est aussi celui d'une vari\'et\'e projective complexe lisse.

 \end{abstract}

\maketitle

\section*{Introduction}

Rappelons une question classique:

\begin{qt}\label{question pb serre}
La classe des groupes kählériens\footnote{Un groupe de présentation finie est dit \emph{kählérien} (resp. \emph{projectif}) s'il peut être réalisé comme le groupe fondamental d'une variété kählérienne compacte (resp. projective lisse).} coïncide-t-elle avec celle des groupes projectifs ?
\end{qt}

Nous d\'emontrons ci-dessous le r\'esultat partiel suivant:

\begin{thm}\label{kvsp} Si $G\subset \mathrm{GL}_n(\Bbb C)$ est l'image d'une repr\'esentation lin\'eaire d'un groupe k\"ahlérien $\pi_1(X)$, il existe une vari\'et\'e projective complexe lisse $X'$ et un sous-groupe d'indice fini de $G$ qui est l'image d'une repr\'esentation lin\'eaire de $\pi_1(X')$. En particulier, un groupe kählérien linéaire est virtuellement projectif.

\end{thm}

La motivation initiale de la question \ref{question pb serre} est le \textsl{probl\`eme de Kodaira} : une vari\'et\'e k\"ahl\'erienne compacte est-elle toujours d\'eformation d'une vari\'et\'e projective complexe lisse? Une reponse affirmative\footnote{Connue en dimensions $1$ et $2$ (Kodaira).} \`a ce probl\`eme impliquerait en particulier qu'une vari\'et\'e k\"ahl\'erienne compacte est toujours diff\'eomorphe \`a une vari\'et\'e projective, et une r\'eponse affirmative \`a la question pr\'ec\'edente.    Les contre-exemples fournis par C. Voisin (\cite{V04} , \cite{V06}) concernent cependant l'aspect bim\'eromorphe du prob\`eme, et conduisent \`a le reformuler seulement pour les mod\`eles minimaux (en g\'en\'eral singuliers) des vari\'et\'es k\"ahl\'eriennes compactes \`a fibr\'e canonique pseudo-effectif\footnote{ces mod\`eles minimaux sont les espaces k\"ahl\'eriens compacts \`a singularit\'es terminales et fibr\'e canonique semi-ample. Lorsque $K_X$ n'est pas pseudo-effectif, $X$ devrait \^etre unir\'egl\'ee et le quotient rationnel ramènerait le calcul du groupe fondamental au cas o\`u $K_X$ est pseudo-effectif.}. Puisque le passage au mod\`ele minimal n'affecte pas le groupe fondamental, cette reformulation du probl\`eme de Kodaira n'affecte pas non plus la conjecture de coïncidence des groupes k\"ahlérien avec leurs analogues projectifs.

\

La d\'emonstration du th\'eor\`eme \ref{kvsp} consiste \`a se ramener, gr\^ace au th\'eor\`eme \ref{rappel th CCE} ci-dessous, \'etabli dans \cite{CCE}, au cas o\`u $X$ est munie d'une submersion holomorphe dont les fibres sont des tores, sur une vari\'et\'e complexe projective lisse $Y$.

\begin{thm}\label{rappel th CCE}
Soit $Z$ une variété kählérienne compacte et $\rho:\pi_1(Z)\to \mathrm{GL}_N(\CC)$ une représentation linéaire de son groupe fondamental. Quitte à remplacer $Z$ par un revêtement étale fini, la variété de Shafarevich $Sh_\rho(Z)$ associée à $\rho$ est (biméromorphe à) une vari\'et\'e k\"ahl\'erienne compacte $X$ qui est l'espace total d'une fibration lisse en tores sur une variété de type général
$$Z\stackrel{sh_\rho}{\To} X:=Sh_\rho(Z)\stackrel{s_\rho}{\To} S_\rho(Z).$$
\end{thm}
Nous renvoyons à \cite{CCE} pour les notions utilisées ci-dessus. La vari\'et\'e $X$ est bim\'eromorphe \`a la vari\'et\'e de Shafarevich de $Z$, obtenue à partir de $Z$ en contractant les sous-vari\'et\'es sur le $\pi_1$ desquelles $\rho$ est triviale. Les tores de la submersion $X\to S_\rho(Z)$ sont les sous-vari\'et\'es maximales sur le $\pi_1$ desquelles $\rho$ a une image abélienne dans $G$.\\

Nous r\'esolvons ensuite dans ce cas tr\`es particulier le probl\`eme de Kodaira. Nous montrons en effet ci-dessous que $X$ admet (apr\`es changement de base fini) une d\'eformation sur sa fibration jacobienne $J(X)$  (proposition \ref{deformation jacobienne relative}). Il est alors possible de montrer que $J(X)$ admet des petites déformations au dessus de la variété projective $Y$ dont les fibres sont des variétés abéliennes (th\'eor\`eme \ref{famille tores}) : ces petites déformations sont donc algébriques.

\begin{merci}
Nous remercions vivement C. Voisin pour son aide dans l'élaboration de certains arguments présents dans ce travail. Le second auteur bénéficie du soutien du CNRS via un séjour de recherche dans l'UMI IMPA/CNRS et tient à remercier l'IMPA pour l'accueil qu'il a reçu.
\end{merci}

\section{Une version relative du critère de Buchdahl}

L'ingrédient principal est une version relative du critère de N. Buchdahl \cite{B06,B08}. Dans les travaux en question, ce critère était utilisé pour donner une démonstration du résultat suivant (dû à Kodaira) ne s'appuyant pas (ou peu) sur la classification des surfaces kählériennes : toute surface kählérienne compacte admet des petites déformations projectives.
\begin{thm}\label{buchdahl relatif}
Soit $\pi:X\To B$ une famille lisse, l'espace total de la famille étant supposé être une variété kählérienne compacte. Considérons une déformation relative de $X/B$ de la forme
$$\xymatrix{\mathcal{X}\ar[rr]^g\ar[rd] & &U\times B\ar[ld]^{\mathrm{pr}_1}\\
&U&
}$$
avec $U$ lisse (au voisinage du point $o\in U$ paramétrant $X$). Soit $\omega$ une métrique kählérienne sur $X$ et supposons que la composition de l'application induite par le produit extérieur et de l'application de Kodaira-Spencer relative
$$T_{U,o}\To \mathrm{H}^1(X,T_{X/B})\stackrel{\wedge [\omega]}{\To} \mathrm{H}^0(B,R^2\pi_*\mathcal{O}_X)$$
est surjective. La famille $X/B$ admet alors des déformations $\pi_u:\mathcal{X}_u\To B$ (avec $u$ arbitrairement proche de $o\in U$) dont toutes les fibres sont projectives. Plus précisément, il existe sur $\mathcal{X}_u$ une classe rationnelle de type (1,1) dont la restriction aux fibres de $\pi_u$ est une classe kählérienne.
\end{thm}
Nous renvoyons à la référence \cite[Chap. 3.4.2]{S06} pour les notions de théories des déformations utilisées ici.\\

Plusieurs remarques s'imposent. Tout d'abord, l'application de produit extérieur ci-dessus est en fait la composition du morphisme naturel
$$\mathrm{H}^1(X,T_{X/B})\To \mathrm{H}^0(B,R^1\pi_*T_{X/B})$$
avec le produit extérieur par la classe $[\omega]$ vue comme un élément de $$[\omega]\in\mathrm{H}^0(B,R^1\pi_*\Omega^1_{X/B}).$$
La condition peut aussi s'exprimer en examinant la composition :
$$T_{U,o}\To \mathrm{H}^1(X,T_X)\stackrel{\wedge \omega}{\To} \mathrm{H}^2(X,\mathcal{O}_X)\To \mathrm{H}^2(X_b,\mathcal{O}_{X_b}),$$
cette dernière devant alors être surjective sur l'espace
$$\mathrm{Im}\left(\mathrm{H}^2(X,\mathcal{O}_X)\To \mathrm{H}^2(X_b,\mathcal{O}_{X_b})\right).$$
D'autre part, d'après la version relative du critère de projectivité de Kodaira (voir \cite{C06}), les fibres de $\pi$ sont d'ores et déjà projectives si l'application de restriction ci-dessus est identiquement nulle.\\
Enfin, il est important de noter que, dans une telle déformation relative, la structure complexe de $B$ reste inchangée et que, d'un point de vue ensembliste, l'application $\pi_u$ n'est autre que l'application $\pi$ (la structure complexe de $X$ n'est modifiée que dans la direction des fibres).\\

\begin{demo}[du théorème \ref{buchdahl relatif}]
Les syst\`emes locaux $R^2 (\mathrm{pr}_1\circ g)_* \QQ_{\mathcal{X}}$ et $R^2 g_* \QQ_{\mathcal{X}}$ sont munis par la construction de \cite[Chap. 10]{V02} d'une structure  de $\QQ$-variation de structure de Hodge, cette structure admettant une polarisation r\'eelle si $g$ en admet une - ce qui n'est pas n\'ecessairement le cas. M\^eme en l'absence d'une polarisation, nous pouvons munir $(\mathrm{pr}_1)_* R^2g_* \QQ_{\mathcal{X}}$ d'une structure de $\QQ$-variation de structure de Hodge telle que l'\'epimorphisme canonique
$$s: R^2 (\mathrm{pr}_1\circ g)_* \QQ_{\mathcal{X}} \to (\mathrm{pr}_1)_* R^2g_* \QQ_{\mathcal{X}}$$
donn\'e par la d\'eg\'erescence de la suite spectrale de Leray pr\'eserve les filtrations de Hodge. La filtration de Hodge $F$ de $(\mathrm{pr}_1)_* R^2g_* \QQ_{\mathcal{X}}$ est juste l'image de la filtration de
Hodge de $R^2 (\mathrm{pr}_1\circ g)_* \QQ_{\mathcal{X}}$. Pour voir que cette filtration $F$ est bien de Hodge (ce qui r\'esulte de $F^i \cap \overline{F^{3-i}} = \{0\}$), il suffit de choisir $b\in B$ et d'observer que la compos\'ee $j\circ s$ o\`u $j$ est le \emph{monomorphisme} naturel de $ (\mathrm{pr}_1)_* R^2g_* \QQ_{\mathcal{X}}$
dans $R^2 h^b_* \QQ_{\mathcal{X}_b}$ (avec $h^b: \mathcal{X}_b=g^{-1}(U\times\{b\}) \to U$) est un morphisme de structure de Hodge.

Comme dans \cite[Prop. 3.3]{V05},  le th\'eorème r\'esulte alors de la proposition 17.20 p. 410 de \cite{V02} appliqu\'ee \`a la $\QQ$-variation de structure de Hodge sous-jacente au syst\`eme local $(\mathrm{pr}_1)_* R^2g_* \QQ_{\mathcal{X}}$.
\end{demo}

\section{Application aux familles en tores}

Dans la suite nous ferons la distinction suivante : une fibration en tores sera une submersion propre dont les fibres sont des tores complexes alors qu'une fibration jacobienne sera une famille en tores munie d'une section (ce qui revient à prescrire l'origine). Nous adopterons également la terminologie suivante.
\begin{defi}\label{defi jac relative}
Soit $X/B$ une famille en tores (avec $X$ une variété kählérienne compacte). Nous appellerons jacobienne relative de $X/B$ la fibration jacobienne\footnote{il n'y a \emph{a priori} pas unicité de la section mais la variété $Y$ est cependant unique.} $Y/B$ ayant les mêmes fibres que celles de $X/B$.
\end{defi}
Notons que la jacobienne relative est une vari\'et\'e k\"ahl\'erienne par un th\'eor\`eme de Varouchas puisque l'on dispose d'un morphisme surjectif propre et lisse $X\times_B  X \to Y$.

La famille $Y/B$ est donc obtenue à partir de $X/B$ en effectuant des translations dans les fibres pour garantir l'existence d'une section holomorphe. Par exemple, si $X$ est un tore de dimension 2 et de dimension algébrique $a(X)=1$ (voir \cite{BL}), la réduction algébrique de $X$ est une famille en tores $X/B$ sur une courbe elliptique $B$ et de fibre $E$ (dont la structure complexe est fixée). La jacobienne relative de $X/B$ n'est autre que $E\times B$.\\

Nous utiliserons à diverses reprises la description suivante de $Y/B$ : si $X/B$ est une famille en tores, considérons $\mathrm{H}_\ZZ$ le système local formé par les $\mathrm{H}_1$ des fibres et $\mathcal{E}$ le sous-fibré holomorphe formé par la partie de type $(-1,0)$ de $\mathrm{H}_\ZZ\otimes \mathcal{O}_B$. La famille $Y/B$ s'obtient naturellement en considérant le quotient de l'espace total de $\mathcal{E}$ par le réseau $\mathrm{H}_\ZZ$ (la section nulle passant au quotient). En particulier, si $\mathcal{Y}$ désigne le faisceau des sections de $Y/B$, il s'insère dans la suite exacte suivante :
\begin{equation}\label{suite exacte jac relative}
\xymatrix{0\ar[r] & \mathrm{H}_\ZZ\ar[r] & \mathcal{E}\ar[r] & \mathcal{Y}\ar[r] & 0}
\end{equation}
Enfin, remarquons qu'une famille en tores $X/B$ définit naturellement une classe dans le groupe $\mathrm{H}^1(B,\mathcal{Y})$. En effet, si $s_i$ et $s_j$ sont deux sections locales de $X/B$ (définies sur des ouverts $U_i$ et $U_j$), alors la différence $s_{i,j}:=s_i-s_j$ est bien définie (indépendamment du choix d'une origine dans les fibres) et se trouve être une section locale de $Y/B$ sur $U_i\cap U_j$. La collection des $s_{i,j}$ est un cocycle de \v{C}ech de degré 1 et fournit bien une classe
$$e_h(X/B)\in \mathrm{H}^1(B,\mathcal{Y}).$$

\subsection{Topologie des familles en tores}
Nous allons maintenant décrire la topologie des familles en tores (avec espace total kählérien) et montrer qu'une famille en tores se déforme vers sa jacobienne relative (après changement de base étale fini). Certains des résultats montrés ci-dessous sont également présent dans le travail de Nakayama \cite{Nak}.\\

Soit donc $\pi:X\To B$ une famille en tores (la variété différentiable sous-jacente aux fibres étant toujours notée $F$) dont l'espace total est une variété kählérienne compacte. La dégénérescence de la suite spectrale de Leray (\emph{via} le théorème des cycles invariants) implique facilement (voir \cite{Cl10}) que la suite des groupes fondamentaux :
\begin{equation}\label{suite exacte gf}
1\To \pi_1(F)\To \pi_1(X)\To\pi_1(B)\To 1
\end{equation}
est exacte. Rappelons-en brièvement la démonstration. Considérons la fibration (de fibre $F$) $\tilde{\pi}:\tilde{X}\To \tilde{B}$ obtenue par image réciproque sur $\tilde{B}$ revêtement universel de $B$. La suite exacte d'homotopie pour cette fibration s'écrit :
$$0\To K\To \pi_1(F)\To \pi_1(\tilde{X})\To 0$$
où $K$ est un sous-groupe de $\pi_1(F)$, en particulier abélien sans torsion. En la dualisant, nous obtenons :
$$0\To \mathrm{H}^1(\tilde{X},\QQ)\To \mathrm{H}^1(F,\QQ)\To \mathrm{Hom}(K,\QQ)\To 0.$$
Or, le théorème des cycles invariants (voir par exemple \cite[th. 16.18, p. 381]{V02} affirme que l'application de restriction $\mathrm{H}^1(\tilde{X},\QQ)\to \mathrm{H}^1(F,\QQ)$ est surjective. Comme $K$ est abélien sans torsion, cela implique que $K$ est trivial et le groupe fondamental de $F$ s'injecte dans celui de $X$ (dans un revêtement l'application est injective au niveau du $\pi_1$).\\ 

Il est bien connu qu'une suite exacte comme (\ref{suite exacte gf}) détermine une classe d'extension
$$e(X/B)\in\mathrm{H}^2(\pi_1(B),\pi_1(F)).$$
D'autre part, si l'on considère la jacobienne relative $Y\To B$ (voir ci-dessus), la famille initiale détermine un élément de $e_h(X/B)\in\mathrm{H}^1(B,\mathcal{Y})$. Ces deux classes de cohomologie sont reliées entre elles de la façon suivante. La suite exacte longue de cohomologie associée à la suite courte (\ref{suite exacte jac relative}) fournit une flèche naturelle
$$\mathrm{H}^1(B,\mathcal{Y})\To \mathrm{H}^2(B,\mathrm{H}_\ZZ).$$
Notons $e_{top}(X/B)$ l'image de $e_h(X/B)$ dans $\mathrm{H}^2(B,\mathrm{H}_\ZZ)$. Sous l'injection naturelle
$$\mathrm{H}^2(\pi_1(B),\pi_1(F))\hookrightarrow \mathrm{H}^2(B,\mathrm{H}_\ZZ),$$
l'image de $e(X/B)$ n'est autre que la classe $e_{top}(X/B)$.\\

Il r\'esulte du th\'eorème principal de \cite{A11} que les classes $e$ et $e_{top}$ sont des classes de torsion (c'est à nouveau une conséquence de la dégénérescence\footnote{dans notre cas, l'argument d'Arapura se réduit bien à cela. En effet, la classe d'extension $e(X/B)$ est (au signe près) $d_2(\mathrm{id})$ où $d_2$ est la différentielle
$$d_2:E^{1,0}_2=\mathrm{H}^0(\pi_1(B),\mathrm{H}^1(\pi_1(F),V))=\mathrm{Hom}_{\pi_1(X)}(V,V)\To \mathrm{H}^2(\pi_1(B),V)=E^{2,0}_2$$
de la suite spectrale de Hochschild-Serre déduite de la suite exacte (\ref{suite exacte gf}) pour le $\pi_1(X)$-module $V:=\pi_1(F)$. Le groupe $\pi_1(F)$ étant abélien, $V$ est en fait un système local défini sur $B$. La suite spectrale de Leray associée au morphisme $f:X\to B$ et au faisceau $f^*V_\QQ$ dégénère en $E_2$ et ceci montre bien que $e(X/B)$ est nulle en cohomologie à coefficients rationnels ou, en d'autres termes, que c'est bien une classe de torsion.} de la suite spectrale de Leray). Cette remarque suffit à montrer que la suite exacte (\ref{suite exacte gf}) est scindée après revêtement étale fini de la base (et même que $X/B$ admet une section différentiable après ce même changement de base). En effet, si $n$ désigne l'ordre de torsion de la classe $e$, la suite longue associée à la suite courte de système locaux
$$0\To \mathrm{H}_\ZZ\stackrel{\times n}{\To} \mathrm{H}_\ZZ \To \mathrm{H}_\ZZ/n\mathrm{H}_\ZZ\To 0$$
montre que la classe $e$ provient en fait de
$$\mathrm{H}^1(B,\mathrm{H}_\ZZ/n\mathrm{H}_\ZZ)=\mathrm{H}^1(\pi_1(B),\mathrm{H}_\ZZ/n\mathrm{H}_\ZZ).$$
Le groupe abélien sous-jacent à $\mathrm{H}_\ZZ/n\mathrm{H}_\ZZ$ étant fini, nous pouvons considérer un sous-groupe d'indice fini $G$ de $\pi_1(B)$ qui agit trivialement sur $\mathrm{H}_\ZZ/n\mathrm{H}_\ZZ$. Le groupe de cohomologie correspondant s'identifie alors à
$$\mathrm{H}^1(G,\mathrm{H}_\ZZ/n\mathrm{H}_\ZZ)=\mathrm{Hom}(G,\mathrm{H}_\ZZ/n\mathrm{H}_\ZZ).$$
En considérant le revêtement correspondant au noyau du morphisme induit par la classe $e_{\vert G}$, on obtient le scindage voulu. Résumons la discussion précédente dans l'énoncé suivant.

\begin{prop}\label{topologie famille tore}
Soit $X/B$ une famille en tores dont l'espace total est un variété kählérienne compacte. Après changement de base étale fini, la classe $e_{top}(X/B)$ s'annule et cette famille possède une section (différentiable) ; en particulier, le groupe fondamental de $X$ devient un produit semi-direct
$$\xymatrix{1\ar[r]& \pi_1(F)\ar[r]& \pi_1(X)\ar[r]&\pi_1(B)\ar@/_/[l]\ar[r]& 1}$$
\end{prop}

\begin{rem}\label{rq non kahler}
L'hypothèse concernant le caractère kählérien de l'espace total est bien entendu essentiel comme le montre l'exemple des surfaces de Hopf (celles admettant une fibration elliptique sur $\PP^1$) ou celui de la variété d'Iwasawa (qui admet une fibration elliptique sur un tore complexe de dimension 2).
\end{rem}

\subsection{Famille en tore et fibration jacobienne}

La classe $e(X/B)$ mesure également l'obstruction à déformer $X/B$ vers sa jacobienne relative $Y/B$.
\begin{prop}\label{deformation jacobienne relative}
Soit $X/B$ une famille de tores. Si $e(X/B)=0$, $X/B$ et $Y/B$ sont déformations l'une de l'autre.
\end{prop}
\begin{demo}
Si la classe $e(X/B)$ s'annule dans $\mathrm{H}^2(B,\mathrm{H}_\ZZ)$, la suite longue associée à la suite (\ref{suite exacte jac relative}) montre que $e_h(X/B)$ provient de $\mathrm{H}^1(B,\mathcal{E})$. Or, une classe dans le groupe $\mathrm{H}^1(B,\mathcal{E})$ correspond à une extension
$$\xymatrix{0\ar[r] & \mathcal{E}\ar[r] & \mathcal{F}\ar[r] & \mathcal{O}_B\ar[r] & 0}$$
du fibré $\mathcal{E}$ par le fibré trivial. A nouveau, nous pouvons faire opérer le réseau $\mathrm{H}_\ZZ$ sur $\mathbb{F}$ l'espace total de $\mathcal{F}$ par translation. La variété complexe obtenue admet naturellement une submersion (propre) sur l'espace total de $\mathcal{O}_B$ :
$$\epsilon:\mathbb{F}/\mathrm{H}_\ZZ\To \CC.$$
Il est clair que le fibre de $\epsilon$ au dessus de 0 n'est autre que $Y$ ; il n'est pas moins évident que celle au dessus de 1 s'identifie à $X$.
\end{demo}

\subsection{Déformations relatives des fibrations jacobiennes}

\begin{thm}\label{famille tores}
Soit $\pi:X\To B$ une fibration jacobienne, $X$ étant toujours supposée kählérienne compacte. La famille $X/B$ admet des petites déformations de telle sorte que les fibres de $\pi$ soient des variétés abéliennes. En particulier, si $B$ est une variété projective, $X$ admet des petites déformations projectives.
\end{thm}
Avant de donner la démonstration de ce dernier résultat, énonçons-en la conséquence suivante qui fût la motivation principale de cette note.
\begin{cor}\label{cor tores}
Si $\pi:X\To B$ est une famille en tores et si $B$ est de plus une variété projective, après revêtement étale fini $X$ peut se déformer en une variété projective ; en particulier, son groupe fondamental est virtuellement celui d'une variété projective.
\end{cor}
\begin{demo}[du corollaire \ref{cor tores}]
Il suffit dans un premier temps d'appliquer la proposition \ref{topologie famille tore} pour annuler $e_{top}(X/B)$ après un revêtement fini de $B$. La proposition \ref{deformation jacobienne relative} montre que dans cette situation $X/B$ se déforme à sa jacobienne relative $Y/B$. Cette famille étant munie d'une section, appliquons le théorème \ref{famille tores} à la famille $Y/B$ : nous pouvons déformer $Y$ en une famille en tores dont les fibres et la base sont algébriques. L'espace total de cette déformation est alors algébriquement connexe (car la fibration est toujours munie d'une section) donc algébrique\footnote{Nous pouvons également raisonner de la manière suivante : 
les arguments ci-dessus montrent qu'après déformation il existe une classe $\alpha\in H^2(X,\QQ)$ dont la restriction aux fibres est celle d'une classe de Hodge. Or, dans le cas des familles de tores, la dégénérescence de la suite de Leray fournit en fait une décomposition canonique en sous-structures de Hodge
$$H^k(X,\QQ)\simeq \bigoplus_{i+j=k}H^i(B,R^jf_*\QQ)$$
(voir \cite{V11}). Nous pouvons donc supposer que $\alpha$ est de type (1,1) et la classe $\alpha+f^*(\omega_B)$ est donc une classe de Hodge (avec $\omega_B$ polarisation sur $B$).} d'après \cite{C81}.
\end{demo}

\begin{demo}[du théorème \ref{famille tores}]
Comme nous considérons une fibration jacobienne avec espace total kählérien, les déformations relatives de la structure complexe sont les mêmes que celles du système local $\mathcal{H}_R:=(\mathrm{H}_1(X_b,\RR))_{b\in B}$. Nous cherchons ici à déformer la structure complexe $I$ donnée par la théorie de Hodge usuelle
\begin{equation}\label{structure usuelle}
(\mathcal{H}^*_\CC)^{1,0}:=\pi_*\Omega^1_{X/B}.
\end{equation}
Nous allons montrer que d'une part certains éléments de
$$A_\RR=\mathrm{H}^0(B,\mathrm{End}(\mathcal{H}_\RR))$$
s'identifient à des déformations infinitésimales et d'autre part que ces déformations ne sont pas obstruées. Il nous restera à identifier leurs images dans $\mathrm{H}^2(X_b,\mathcal{O}_{X_b})$ pour constater que le critère \ref{buchdahl relatif} s'applique.

Soit donc $\varphi$ un élément de $A_\RR$. Nous souhaitons déformer la structure complexe $I$, c'est-à-dire construire un développement (au moins formel) de la forme :
$$I_t=I+t\varphi_1+t^2\varphi_2+\dots$$
avec $I_t^2=-1$ et $\varphi_1=\varphi$. La condition à l'ordre 1 est la suivante : $\varphi$ doit anti-commuter avec $I$. De façon plus intrinsèque, le théorème des cycles invariants permet de munir $A_\RR$ d'une structure de Hodge de poids 0 et les éléments de $A_\RR\cap(A_\CC^{1,-1}\oplus A_\CC^{-1,1})$ sont exactement ceux qui anti-commutent avec $I$. Il nous reste à vérifier que ces déformations infinitésimales ne sont pas obstruées (ceci est évident dans le cas des tores complexes mais ici nous devons prendre en compte l'action du groupe de monodromie) ; c'est l'objet du lemme \ref{structure complexe syst local} ci-dessous.

Les déformations infinitésimales comme ci-dessus peuvent s'obtenir de la façon suivante ; considérons pour cela l'application de restriction
$$j_b^*:\mathrm{H}^2(X,\RR)\To \mathrm{H}^2(X_b,\RR)$$
et remarquons que, $X_b$ étant un tore, ce dernier groupe de cohomologie s'identifie à
$$\mathrm{H}^2(X_b,\RR)=\bigwedge^2 \mathrm{H}_1(X_b,\RR)^*.$$
Nous pouvons donc considérer les éléments de l'image de $j_b^*$ comme des formes bilinéaires sur le système local $\mathcal{H}_\RR$ ou encore comme des morphismes de $\mathcal{H}_\RR\To\mathcal{H}^*_\RR$. Enfin, le choix d'une classe kählérienne globale $[\omega]$ fournit un isomorphisme de $\mathcal{H}_\RR$ vers son dual et nous pouvons considérer des éléments de la forme
$$\varphi=[\omega]^{-1}\circ\eta$$
avec $\eta\in \mathrm{Im}(j_b^*)$. Si la partie de type (0,2) de l'image de $j_b^*$ est nulle, nous savons déjà que les fibres sont projectives (voir \cite{C06}). Dans le cas contraire, la structure complexe du système local admet des déformations et il est immédiat de constater que la contraction de $\varphi$ avec la classe $[\omega]$ n'est autre que $\eta^{0,2}$. Le théorème \ref{buchdahl relatif} s'applique donc et $X$ admet des déformations (arbitrairement petites) au dessus de $B$ dont les fibres sont projectives.
\end{demo}

Nous avons utilisé ci-dessus le lemme suivant.
\begin{lem}\label{structure complexe syst local}
Reprenons les notations ci-dessus. Les déformations infinitésimales données par $A_\RR\cap(A_\CC^{1,-1}\oplus A_\CC^{-1,1})$ ne sont pas obstruées.
\end{lem}
\begin{demo}
Considérons le groupe $G:=\mathrm{Aut}(\mathcal{H}_{\RR})$. Le groupe $G$ est le groupe des inversibles de l'algèbre $A_\RR:=\mathrm{End}(\mathcal{H}_{\RR})$ et son complexifié $G_{\CC}$ est le groupe des endomorphismes inversibles de $A_\CC:=\mathrm{End}(\mathcal{H}_{\CC})$. Pour $g\in G$ on peut construire $I^g=gIg^{-1}$ qui est une structure complexe sur $\mathcal{H}_\RR$ définissant automatiquement une famille holomorphe (\emph{i.e.} la partie de type $(1,0)$ est un sous-fibré holomorphe de $\mathcal{H}_{\CC}$) et $X^g \to B$ est une famille de tores sur $B$.

En effet, la connexion plate sur $\mathcal{H}_{\RR}$ définit un scindage diff\'erentiable $T_X= \pi^* T_B \oplus T_{X|B}\simeq \pi^*T_B \oplus \pi^* \mathcal{H}_{\RR}$. La structure presque complexe de $X$ se scinde de fa\c{c}on compatible sous la forme $\pi^* J_B \oplus \pi^* J_{X|B}$
o\`u $J_B$ est la structure presque complexe de $B$ et $J_{X|B}$ est la structure complexe (\ref{structure usuelle}) sur $\mathcal{H}_{\RR}$.
Le revêtement universel $\tilde X \to X$ est (par transport parall\`ele bas\'e en $\tilde b \in B$)  diff\'eomorphe \`a $\tilde B \times \mathcal{H}_{\RR, \tilde b}$ et sa structure presque complexe se scinde de fa\c{c}on compatible en $J_{\tilde X}= J_{\tilde B}\oplus J_{\mathcal H}$ o\`u les deux
structures presque complexes ne d\'ependent que des param\`etres de $\tilde B$.

Si $g\in \mathrm{Aut}(\mathcal{H}_{\RR})\subset \mathrm{Aut}(\mathcal{H}_{\RR, \tilde b})$ on d\'efinit une structure presque complexe $J^g$ sur $X$ par la formule $J^g:=\pi^* J_B \oplus \pi^* g. J_{X|B}.g^{-1}$. Son rel\`evement $J^g_{\tilde X}$  \`a $\tilde X$ est manifestement \'egal \`a $\phi^{-1}. J_{\tilde X}$ o\`u $\phi: \tilde X \to \tilde X$ est
le diff\'eomorphisme $(b, h)\mapsto (b, g.h)$ - on utilise ici de fa\c{c}on cruciale que $g$ ne d\'epend pas de $b$. Le th\'eor\`eme de Newlander-Nirenberg implique alors que $J^g_{\tilde X}$ donc $J^g$ est int\'egrable. 

Les endomorphismes de $\mathcal{H}_{\RR}$ commutant \`a $I$ forment exactement la partie de type $(0,0)$ de la structure de Hodge réelle de poids $0$ sur $A_\RR$ (\textsc{sh} dont les composants possiblement non triviaux sont de type (-1,1), (0,0) et (1,-1)). Il s'ensuit que le groupe $G^0$ des automorphismes commutant \`a $I$ n'est autre que $G\cap F^0 A_{\CC}$ qui est un ouvert dans $A_{\RR}^{0,0}$. Le groupe des inversibles de $F^0 A_{\CC}$ est lui noté $F^0 G_{\CC}$.

On constate alors que $G/G^0 \To G_{\CC} / F^0 G_{\CC}$ est un homéomorphisme local ce qui fournit une structure complexe sur $D=G/G^0$. De plus, l'espace tangent (au point base correspondant à $I$) à $D$ s'identifie bien aux déformations infinitésimales $A_\RR\cap(A_\CC^{1,-1}\oplus A_\CC^{-1,1})$.

On regarde l'espace total $\mathrm{H}_{\CC}  / B$ du complexifié de $\mathcal{H}_{\RR}$. Les translations par le sous-fibré holomorphe $V$ et le système local $\mathrm{H}_{\ZZ}$ fournissent une relation d'équivalence holomorphe (qui localement sur $B$ est la relation d'équivalence sous-jacente à une action de groupe) dont le quotient est $X/B$. On définit alors un automorphisme holomorphe $\phi$ de fibrés vectoriels $(\mathcal{H}_{\CC}/B) \times G_{\CC}$  par $(h,g)\mapsto (gh,h)$ ce qui permet de définir le sous-fibré holomorphe $\phi\cdot V \subset \mathrm{H}_  {\CC}\times G_\CC / B\times G_\CC$. On peut faire agir le système local $\mathrm{H}_{\ZZ}$ sur $\mathcal{H}_{\CC}\times G_\CC/ B\times G_\CC$ sans tordre par $G$ c'est-à-dire en prenant pour chaque $g\in G_\CC$ l'action originelle sur $\mathcal{H}_{\CC}$.  Le quotient de $\mathrm{H}_{\CC}$ par la relation d'équivalence  engendrée par $\phi\cdot V $ et $\mathrm{H}_{\ZZ}$ fournit au dessus d'un voisinage $U$ de l'identité dans $G$ une famille holomorphe de d\'eformations de $X/B$ qui descend manifestement \`a un voisinage de la classe de l'identit\'e dans $D$ et r\'ealise les d\'eformations infinit\'esimales prescrites.
\end{demo}

\begin{qt}
Nous ne savons pas si la conclusion du théorème \ref{cor tores} persiste sans remplacer $X$ par un revêtement étale fini. Comme ce revêtement provient de la base de la fibration $f:X\To B$, il faudrait étudier le problème de façon équivariante, c'est-à-dire lorsqu'un groupe fini $\Gamma$ agit à la fois sur $X$ et $B$ et que $f$ est équivariante pour les actions considérées. L'adaptation à ce cadre équivariant du théorème \ref{famille tores} ne semble en effet pas immédiate.

Notons que le cas où $B$ est un point est traité dans l'article \cite{BR} et montre en particulier que la question \ref{question pb serre} a une réponse positive dans le cas virtuellement abélien : un groupe kählérien virtuellement abélien est aussi le groupe fondamental d'une variété projective lisse.
\end{qt}


\providecommand{\bysame}{\leavevmode\hbox to3em{\hrulefill}\thinspace}

\end{document}